\documentclass[preprint,authoryear]{elsarticle}

\usepackage[utf8]{inputenc}
\usepackage[english]{babel}
\usepackage{amsmath}
\usepackage{amsthm}
\usepackage{xcolor}
\usepackage{amssymb}
\usepackage{hyperref}
\usepackage{comment}

\newcommand{\red}[1]{\textcolor{black}{#1}}

\usepackage{mathtools}
\usepackage{subcaption}
\usepackage{graphics} 

\usepackage{algorithm}
\usepackage{algpseudocode}

\usepackage{wrapfig}

\usepackage{multirow}

\journal{Journal of Computational and Applied Mathematics}

\begin{document}

\begin{frontmatter}
\begin{abstract}
Recently, there has been an increasing
interest in modelling and computation of physical
systems with neural networks. Hamiltonian systems are an elegant and compact formalism in classical
mechanics, where the dynamics is fully determined by one scalar function,
the Hamiltonian. The solution trajectories are often constrained to evolve on a submanifold of a linear vector space. In this work, we propose new approaches for the accurate approximation of the Hamiltonian function of constrained mechanical systems given sample data information of their solutions. We focus on the importance of the preservation of the constraints in the learning strategy by using both explicit Lie
group integrators and other classical schemes.
\end{abstract}
\title{Learning Hamiltonians of constrained mechanical systems}
\author[1]{Elena Celledoni}
\author[1]{Andrea Leone}
\author[1]{Davide Murari}
\ead{davide.murari@ntnu.no}
\author[1]{Brynjulf Owren}

\address[1]{Dept. of Mathematical Sciences, Norwegian University of Science and Technology}




\end{frontmatter}
\section{Introduction}\label{se:intro}
Neural networks have been proven to be effective in learning patterns from data in many different contexts. Recently there has been an increasing interest in applying neural networks to learn physical models from data, for example models of classical mechanics. For Hamiltonian systems, multiple approaches have been proposed to approximate the energy function, see, e.g., \cite{chen2019symplectic}, \cite{greydanus2019hamiltonian},  \cite{zhong2019symplectic}, \cite{finzi2020simplifying}, \cite{offen2021symplectic}. Building on these results, we propose an improved learning procedure. Our main contribution is an approach to learn the Hamiltonian for systems defined on \red{the cotangent bundle $T^*\mathcal{Q}$ of some manifold $\mathcal{Q}$ embedded in a vector space. Under the assumption that $T^*\mathcal{Q}$ is homogeneous, we show how to do that while preserving the geometry during the learning phase.}
In this paper, by preservation of the geometry we mean the accurate conservation of the constraints rather than of other geometric features such as symplecticity, energy or other first integrals of the system.

As in \cite{finzi2020simplifying}, we express the dynamics of constrained systems by embedding the problem in a vector space of larger dimension, but in our approach we do not make use of Lagrange multipliers. With the aim of understanding the importance of the geometry in this approximation problem, we compare learning procedures based on numerical integrators that preserve the phase space of the system with others that do not. We restrict to homogeneous spaces where Lie group methods can preserve the geometry up to machine accuracy (see, e.g., \cite{celledoni2021}). For example, multi-body lumped mass systems fall naturally in this setting \cite[Chapter 2]{lee2017global}. \red{This restriction still includes systems with the configuration manifold} that is a Lie group, as in some problems of rigid body and rod dynamics, but we will not consider these applications here. The experiments show that there are specific problems where approximating the Hamiltonian using a Lie group method can be relevant. Surprisingly, in many other settings classical Runge--Kutta integrators produce comparable results.

The main focus of the present paper is to learn an approximation of a Hamiltonian system where the training data are given as a set of trajectory segments. To do so, one could learn the dynamics either by approximating the Hamiltonian vector field or the Hamiltonian function as done in our work. Another relevant difference in the learning framework consists of considering in the training procedure either one time step of the flow map (see, e.g., \cite{greydanus2019hamiltonian}) or a sequence of successive time steps as proposed in \cite{chen2019symplectic}. In the latter work it is shown with experimental evidence that taking into account temporal dependencies improves performance. We follow the second strategy when dealing with unconstrained systems, whereas we test both of them with our approach to constrained systems.

In principle, the Hamiltonian can be any differentiable function. However, for mechanical systems, it is often made by the sum of (quadratic) kinetic energy and a potential energy, \cite{whittaker04ato}, \cite{hairer2006geometric}, \cite{marsden1995introduction}. Following \cite{zhong2019symplectic}, we make the ansatz that the kinetic energy is characterized by a symmetric and positive definite matrix, and hence we aim to estimate it.

We conclude this Section with a more precise definition of the problem of interest. In the second Section, we introduce the Hamiltonian formalism for both unconstrained and constrained systems. In the third Section, we focus on unconstrained systems, presenting the general learning procedure that will be extended to constrained systems in the fourth Section.  We also discuss how additional known information about the dynamical system can be included in the network training procedure. The experimental results show that physics-based regularization could be helpful to improve the extrapolation capability of the network and its stability in the presence of noise. In the last Section, we formalize the problem of learning a constrained Hamiltonian mechanical system and discuss the importance of the geometry for this class of problems. Finally, we complete this Section with numerical experiments in the \texttt{PyTorch} framework, showing how the predicted Hamiltonian depends on some training parameters \red{and on the presence of noise}. \red{The numerical implementations are available in the GitHub repository associated to the paper\footnote{\url{https://github.com/davidemurari/learningConstrainedHamiltonians}}}.

\subsection{Description of the problem} \label{se:description}
Suppose to be given a set of $N$ sampled trajectories coming from a Hamiltonian system defined on a submanifold \red{$\mathcal{M}=T^*\mathcal{Q}$ of $\mathbb{R}^{2n}$, where $T^*\mathcal{Q}$ is the cotangent bundle of the configuration manifold $\mathcal{Q}$ (see \cite{lee2012smooth}[Chapter 11] for more details).} Moreover, assume that each of these trajectories contains $M$ equispaced (in time) points. In other words, suppose that
\begin{equation}\label{eq:sampledTraj}
\{(x_i,\bar{y}_i^2,...,\bar{y}_i^M)\}_{i=1,...,N},\;\bar{y}_i^j = \Phi^{(j-1)\Delta t}_{X_H}(x_i)
\end{equation}
as a training set, where $\Phi_{X_H}^t$ is the \red{time t-}flow of the exact, unknown Hamiltonian system. In practice, we never have access to the exact trajectories but to either a noisy version of them or a numerical approximation.

The approach we use aims to approximate the vector field $X_H\in\mathfrak{X}(\mathcal{M})$ that governs the dynamics, \red{where by $\mathfrak{X}(\mathcal{M})$ we denote the collection of all smooth vector fields on $\mathcal{M}$.}. However, we know that such a vector field is Hamiltonian, i.e. there exists a scalar function $H:\mathcal{M}\rightarrow\mathbb{R}$ which, together with the geometry given by $\mathcal{M}$, characterizes the dynamics completely. For this reason, we do not need to directly approximate $X_H$, but just $H$ and then eventually recover $X_H$.

The problem under consideration can be described as an inverse problem, since we want to infer the function $H$ from trajectory data of the corresponding dynamical system rather than from samples of the function $H$ itself. This description of the problem motivates how we measure the accuracy of our approximation, denoted by a parametric model $H_{\Theta}$. Indeed, the target is not to approximate the trajectories of the given Hamiltonian system with some neural network, but to approximate the Hamiltonian. Thus the quality of the approximation can be computed in at least two ways.  First, one can compare some measured trajectories with those obtained from the approximation. More precisely, we randomly generate $\tilde{N}$ initial conditions $z_i\in \mathcal{M}$, \red{their $\tilde{M}$ time updates,} and compute
\begin{equation}\label{eq:MSETraj}
\mathcal{E}_1\left(\left\{u_i^j\right\}_{i=1,...,\tilde{N}}^{j=1,...,\tilde{M}}\, ,\, \left\{v_i^j\right\}_{i=1,...,\tilde{N}}^{j=1,...,\tilde{M}}\right) = \frac{1}{\tilde{N}\tilde{M}}\sum_{j=1}^{\tilde{M}}\sum_{i=1}^{\tilde{N}} \left\|u_i^j - v_i^j\right\|^2,
\end{equation}
where $\|\cdot\|$ is the Euclidean norm of $\mathbb{R}^{2n}$,  \red{$u_i^1 = z_i$, $v_i^1=z_i$}, $u_i^{j+1} = \Psi_{X_H}^h(u_i^j)$ and $v_i^{j+1} = \Psi_{X_{H_{\Theta}}}^h(v_i^j)$ for a numerical integrator $\Psi^h$ of choice. \red{One can randomly generate these initial conditions for academic examples where the true Hamiltonian is actually known. In this case, $\tilde{N}$ and $\tilde{M}$ can be specified arbitrarily, usually with $\tilde{N}$ less the the number of training trajectories $N$. On the other hand, in more realistic applications one has to work with the initial conditions for which the related trajectory segments are known. In this case $\tilde{N}$ and $\tilde{M}$ are constrained by the available data, in particular the number of total trajectories is split into $N$ for training and $\tilde{N}$ for test.} In our experiments, we adopted the $\texttt{SciPy}$ implementation of the Dormand-Prince pair of order (5,4) with a strict tolerance. In fact, following the PyTorch implementation of the mean squared error, $\mathcal{E}_1$ is actually divided by $2n$. Alternatively, as introduced in \cite{david2021symplectic}, one can compare pointwise values of the approximated and the true Hamiltonian, when known. This gives 
\begin{equation}\label{eq:avgError}
\mathcal{E}_2(H,H_{\Theta})=\frac{1}{\tilde{N}} \sum_{i=1}^{\tilde{N}} \left | H(z_i) - H_{\Theta}(z_i)-\frac{1}{\tilde{N}}\sum_{l=1}^{\tilde{N}}\left(H(z_l) - H_{\Theta}(z_l)\right)\right|,
\end{equation}
where $\mathcal{E}_2$ handles the fact that Hamiltonians differing, on $\mathcal{M}$, by a constant generate the same vector field. Indeed, $\mathcal{E}_2(H,H+c)=0$.
\section{Hamiltonian mechanical systems}\label{se:ham}
In this work, we focus on Hamiltonian mechanical systems based on a configuration manifold $\mathcal{Q}\subseteq \mathbb{R}^n$. We now introduce some basic elements of the theory of unconstrained Hamiltonian dynamics on $\mathbb{R}^{2n}$, which corresponds to the case $\mathcal{Q}=\mathbb{R}^n$. Then we extend this formulation to constrained systems on $T^*\mathcal{Q}\subset\mathbb{R}^{2n}$.  

The Hamiltonian formalism gives a particular class of conservative vector fields which, in contrast to the Lagrangian one, can always be expressed with a system of first-order ordinary differential equations. For the unconstrained case, the equations are of the form $\dot{x}(t) = \mathbb{J}\nabla H(x(t)):=X_{H}(x(t))$ where $x(t)=[q(t),p(t)]\in\mathbb{R}^{2n}$ comprises the configuration variables and their conjugate momenta. Here, $H:\mathbb{R}^{2n}\rightarrow\mathbb{R}$ is a smooth function called the Hamiltonian of the system, and $\mathbb{J}\in\mathbb{R}^{2n\times 2n}$ is the symplectic matrix.

In this work, we focus on Hamiltonian systems whose energy function is of the form
\[
H(q,p) = \frac{1}{2}p^TM^{-1}(q)p + V(q)
\]
where $M(q)$ is the mass matrix of the system, possibly depending on the configuration $q\in\mathbb{R}^n$, and $V(q)$ is the potential energy of the system. This is not a too restrictive assumption since it still includes a quite broad family of systems. For unconstrained systems, we will further restrict to the case where $M$ is a constant matrix and the Hamiltonian is separable. 
This assumption allows to implement symplectic numerical integration without needing implicit updates. On the other hand, in the constrained setting we aim at preserving the geometry of the numerical flow map rather than other properties such as symplecticity. As a consequence, we can work with variable mass matrices still using explicit numerical integrators as in the unconstrained case.

We now briefly formalize how to extend this formulation to Hamiltonian systems that are holonomically constrained on some configuration manifold $\mathcal{Q}=\{q\in\mathbb{R}^n:\,g(q)=0\}$ embedded in $\mathbb{R}^n$ (for a more detailed derivation of this formalism we refer to \cite[Chapter 8]{lee2017global}). Many mechanical systems relevant for applications are characterized by the presence of some constraints that are coupled to the ODE defining the dynamics.
One way to model this kind of problems is based on Lagrange multipliers, which lead to differential algebraic equations (DAEs).
There has been some work in the direction of extending the Hamiltonian neural network's framework to constrained systems (see, e.g., \cite{finzi2020simplifying} in which this strategy of introducing Lagrange multipliers is applied).

In this manuscript, we want to present an alternative approach based on the assumption that the constrained manifold $\mathcal{Q}$ is embedded in some linear space $\mathbb{R}^n$. This is actually not a restriction, since Whitney's embedding theorem always guarantees the existence of such an ambient space (see, e.g., \cite[Chapter 6]{lee2012smooth}). More explicitly, because of this embedding property, constrained multi-body systems can be modelled by means of some projection operator and the vector field is written in such a way that it directly respects the constraints, without the addition of algebraic equations.

Furthermore, we assume that the components $g_i:\mathbb{R}^n\rightarrow \mathbb{R}$, $i=1,...,m$, are \red{functionally independent on the zero level set}, so that the Hamiltonian is defined on the $(2n-2m)$ dimensional cotangent bundle $\mathcal{M}=T^*\mathcal{Q}$. Working with elements of \red{the tangent space at $q$,} $T_q\mathcal{Q}$\red{,} as vectors in $\mathbb{R}^n$, we introduce a linear operator that defines the orthogonal projection of an arbitrary vector $v\in\mathbb{R}^n$ onto $T_q \mathcal{Q}$, i.e.
\[
\forall q\in \mathcal{Q},\text{ we set }P(q):\mathbb{R}^n\rightarrow T_q\mathcal{Q},\;\;v\mapsto P(q)v.
\]
$P(q)^T$ can be seen as a map sending vectors of $\mathbb{R}^n$ into covectors in $T_q^*\mathcal{Q}$. If $g(q)$ is differentiable, assuming $G(q)$ is the Jacobian matrix of $g(q)$, we have $T_q \mathcal{Q} = \mathrm{Ker}\,G(q)$, and so $P(q) = I_n - G(q)\left(G(q)^TG(q))\right)^{-1}G(q)^T$, where $I_n\in\mathbb{R}^{n\times n}$ is the identity matrix. This projection map allows us to define Hamilton's equations as follows
\begin{equation}\label{eq:ProjectionEq}
\begin{cases}
\dot{q} = P(q)\partial_pH(q,p)\\
\dot{p} = -P(q)^T\partial_qH(q,p) + W(q,p)\partial_pH(q,p),
\end{cases}
\end{equation}
where
\red{
\[
\begin{split}
W(q,p)&=P(q)^T\Lambda(q,p)^T P(q) + \Lambda(q,p)P(q) -P(q)^T\Lambda(q,p)^T,\\ &\text{with}\quad \Lambda(q,p) = \frac{\partial P(q)^Tp}{\partial q}.
\end{split}
\]
}
It is important to remark that since $T^*\mathcal{Q}\subset\mathbb{R}^{2n}$, we can work with the coordinates of the ambient space in the subsequent development. We notice that when $\mathcal{Q}=\mathbb{R}^n$, we can set $P(q)=I$ and recover the unconstrained formulation. These equations of motion can be derived by the standard Hamilton's variational principle on the phase space or by the Legendre transform applied to the Euler-Lagrange equations. However, due to the geometry of the system, the variations need to be constrained to the right spaces and this is done with the projection map $P(q)$. We will focus on the case $Q = S^2\times ... \times S^2 = (S^2)^k$ in Section \ref{se:spherical}, where the mass matrix $M(q)$ and equation (\ref{eq:ProjectionEq}) takes a structured form\red{, with $S^2$ the unit sphere in $\mathbb{R}^3$.}
\section{Learning unconstrained systems}\label{se:hamNet}
\begin{figure}[h!]
    \centering
    \includegraphics[width=\textwidth]{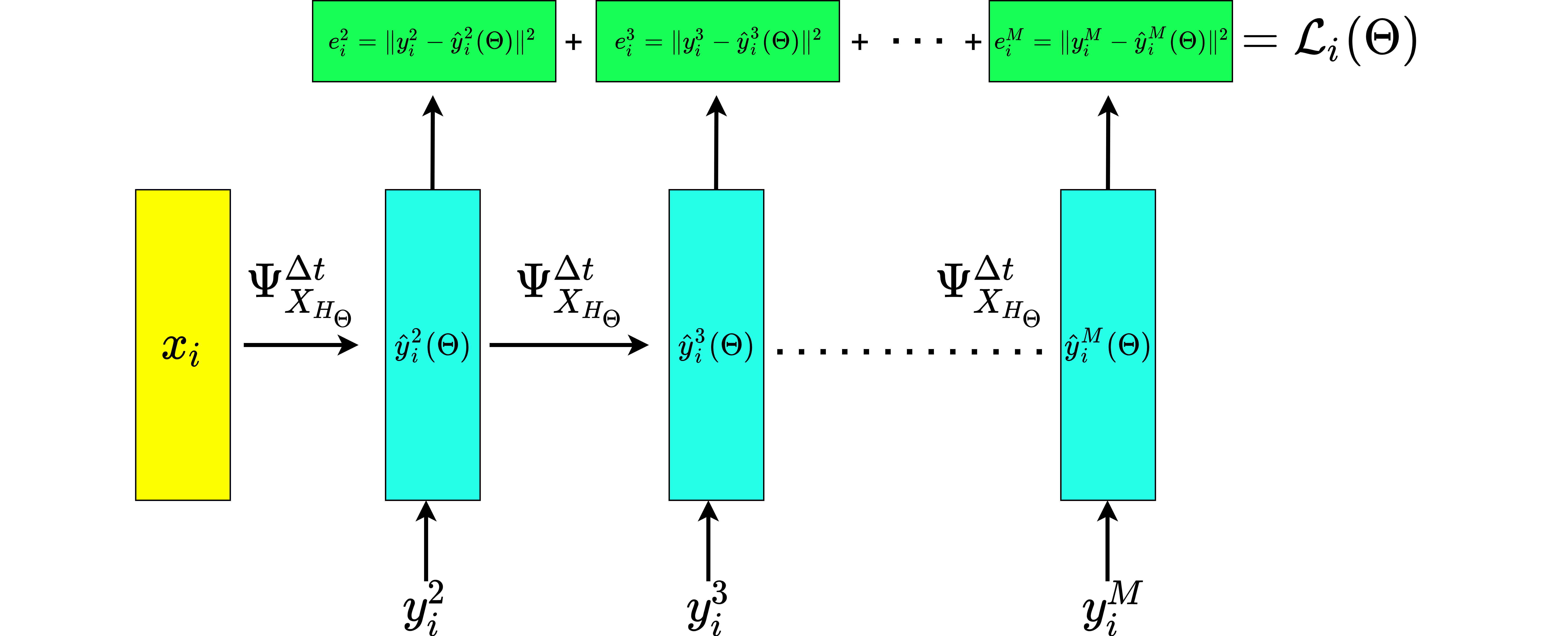}
    \caption{\red{Forward pass of an input training trajectory $(x_i,y_i^2,...,y_i^M)$. The picture highlights the resemblance to an unrolled version of a Recurrent Neural Network. The network outputs $(\hat{y}_i^2,…,\hat{y}_i^M)$.}}
    \label{fig:rhnn}
\end{figure}

As in \cite{chen2019symplectic}, we base the training on a recurrent approach, that is graphically described in Figure \ref{fig:rhnn}.

As mentioned in Subsection \ref{se:description}, we work with numerically generated training trajectories that we denote by \[\{(x_i,y_i^2,...,y_i^M)\}_{i=1,...,N}.\]
We limit the treatment of noisy training data to Subsection \ref{se:noise}.
To obtain an approximation of the Hamiltonian $H$,
we define a parametric model $H_{\Theta}$ and look for a $\Theta$ so that the trajectories generated by $H_{\Theta}$ resemble the given ones. \red{$H_{\Theta}$ in principle can be any parametric function depending on the parameters $\Theta$. In  our approach, $\Theta$ will collect a factor of the mass matrix and the weights of a neural network, as specified in equation \eqref{thetaparam}.} We use some numerical one-step method $\Psi_{X_{H_{\Theta}}}^{\Delta t}$ to generate the trajectories
\begin{equation}\label{eq:psi}
\hat{y}_i^j(\Theta) :=\Psi_{X_{H_{\Theta}}}^{\Delta t}(\hat{y}_i^{j-1}(\Theta)),\quad \hat{y}_i^1(\Theta) := x_i, \quad j=2,\dots,M, \; i=1,\dots,N.
\end{equation}
For unconstrained problems we use symplectic numerical integrators, since they can take an explict form and \red{their adoption in the training procedure allows to have a target modified Hamiltonian to approximate} (see, e.g., \cite{zhu2020deep}).
We then optimize a loss function measuring the distance between the given trajectories $y^j_i$ and the generated ones $\hat{y}_i^j$, defined as
\begin{equation}\label{eq:loss}
\mathcal{L}(\Theta):=\red{\frac{1}{2n}\frac{1}{NM}\sum_{i=1}^N\mathcal{L}_i(\Theta) = }\frac{1}{2n}\frac{1}{NM}\sum_{i=1}^N\sum_{j=1}^M \|\hat{y}_i^j(\Theta)- y_i^j\|^2,
\end{equation}
where $\|\cdot\|$ is the Euclidean metric of $\mathbb{R}^{2n}$. This is implemented with the PyTorch \texttt{MSELoss} loss function. Such a training procedure resembles the one of Recurrent Neural Networks (RNNs), introduced in \cite{rumelhart1985learning}\red{, as shown for the forward pass of a single training trajectory in Figure \ref{fig:rhnn}.}
Indeed, the weight sharing principle of RNNs is reproduced by the time steps in the numerical integrator which are all based on the same approximation of the Hamiltonian, and hence on the same weights $\Theta$. \red{Finally, in Algorithm \ref{alg:one} we report one training epoch for a batch of data points.}

\begin{algorithm}[h!]
    \begin{algorithmic}[1]
        \caption{One epoch of the recurrent approximation of the Hamiltonian.}
        \label{alg:one}
                \State {Choose a numerical integrator ($s$ stages)}
                \State {$\hat{N}\gets \texttt{batch size}$,\qquad $\text{Loss} \gets 0$}
                \For {$i=1,\dots,\hat{N}$}
                \State {$\hat{y}_i^1 \gets x_i$}
                    \For {$j = 1,\dots,M$}
                        \State {$\hat{y}_i^{j,[1]} \gets \hat{y}_i^j$}
                        \For{$k = 1,\dots,s-1$}
                            \State{Compute current value of Hamiltonian $H_{\Theta}(\hat{y}_i^{j,[k]})$}
                            \State Compute $\nabla H_{\Theta}(\hat{y}_i^{j,[k]})$ \Comment{With automatic differentiation}
                            \State{Compute stage $\hat{y}_i^{j,[k+1]}$}
                            \EndFor
                    \State{Compute $\hat{y}_i^{j+1}$}
                    \State{Increase Loss following equation \eqref{eq:loss}}
                    \EndFor
            \EndFor
                \State{Optimize Loss}
            
    \end{algorithmic}
\end{algorithm}
\subsection{Architecture of the network}
In this work, the role of the neural network is to model the Hamiltonian, i.e. a scalar function defined on the phase space $\mathbb{R}^{2n}$. Thus, the starting and arrival spaces are fixed. For unconstrained systems we assume that \[H(q,p) = \frac{1}{2}p^TM^{-1}p + V(q)=K(p)+V(q)\] is separable. Here, the kinetic energy is a quadratic form defined by the symmetric positive definite matrix $M^{-1}$. It can hence be modelled through a learnable matrix $A$\red{, $K(p) \approx K_A(p)$,} by replacing $M^{-1}$ or $M$ with $A^TA$ during the learning procedure. This modelling choice improves extrapolation properties since it allows to \red{learn local (on a compact set) information} that is valid on a larger domain, i.e. the mass matrix. In Section \ref{se:constr} we extend this reasoning to some configuration dependent mass matrices, where $M(q)$ is modelled through a constant symmetric and positive definite matrix. Recalling that $A^TA$ can even be singular or close to singular, one can promote the positive definiteness of the modelled matrix adding a positive definite perturbation matrix to $A^TA$. \red{Notice that, in principle, the imposition of the positive (semi)definiteness of the matrix defining the kinetic energy is not necessary, but it allows to get more interpretable results. Indeed, it is known that the kinetic energy should define a metric on $\mathbb{R}^n$ and the assumption we are making guarantees such a property.} For constrained systems we proceed in a similar way, as shown in equation \eqref{eq:massmat}.  For the potential energy, a possible modelling strategy is to work with standard feedforward neural networks, and hence to define
\[ V(q) \approx V_{\theta}(q) = f_{\theta_m}\circ ...\circ f_{\theta_1}(q),\] \[\theta_i = (W_i,b_i)\in\mathbb{R}^{n_i\times n_{i-1}}\times \mathbb{R}^{n_i},\;\theta:=[\theta_1,...,\theta_m], \]\[ f_{\theta_i}(u) := \Sigma(W_iu + b_i),\;\mathbb{R}^n\ni z\mapsto \Sigma(z) = [\sigma(z_1),...,\sigma(z_n)]\in\mathbb{R}^n, \] for example with $\sigma(x) = \tanh(x)$. In particular applications, where some additional information is known about the system, one can impose more structure on the architecture modelling $V(q)$. For example, in the case of odd potential or rotationally symmetric potential, one can define respectively an odd neural network $V_{\theta}$ or a rotationally equivariant one \red{(see, e.g., \cite{celledoni2021equiv}). Therefore, we have that
\begin{equation} \label{thetaparam}
\Theta = [A, \theta], \quad H(q,p) \approx H_{\Theta}(q,p) = K_A(p) + V_{\theta}(q).
\end{equation}
}
We remark that the Hamiltonian does not need to be approximated by a neural network, and hence in a compositional way. Many other parametrizations are possible. For example, starting from the sparse identification of dynamical systems approach presented in \cite{brunton2016discovering}, in \cite{dipietro2020sparse} it is proposed to parametrize the Hamiltonian with a dictionary of functions, for example polynomials and trigonometric functions. In our work, however, we opt for standard feedforward neural networks as the modelling assumption.

We now provide further details on the extrapolation capabilities of this network model. The learning procedure presented above is based on extracting temporal information coming from a set of trajectories belonging to a compact subset $\Omega\subset \mathbb{R}^{2n}$. In general, there is no reason why the Hamiltonian should be accurate outside of this set. \red{To be more precise, denoting by $T>0$ the largest time at which we know the trajectories}, we have that
\begin{enumerate}
    \item \red{given enough samples in set $\Omega$, distributed in order to capture the behaviour of the dynamical system,} the prediction of the network is expected to be accurate in \red{$\Omega_{[0,T]}:=\{\Phi^t_{X_H}(x):\,t\in [0,T],\,x\in\Omega\}$}, i.e. for any \red{$z_0\in \Omega_{[0,T]}$} and any $\bar{t}>0$ such that \red{$\Phi^t(z_0)\in \Omega_{[0,T]}$} for all $t\in [0,\bar{t}]$,
    \item outside $\Omega_{[0,T]}$ one cannot guarantee that the prediction will be accurate.
\end{enumerate}
If we think of classical regression problems or even classification ones, it seems reasonable not to have information about the approximated quantity outside the sampled area. In those cases, with generalization we mean being sufficiently accurate close to the training points but still inside the sampled domain. However, here we know that the inferred function $H(q,p)$ has physical meaning and properties, so we might incorporate global known information about it to extend the applicability of the predictions.

This discussion supports the architectural choice for the kinetic energy suggested before (as in \cite{zhong2019symplectic}). \red{Indeed, supposing the Hamiltonian is separable, we know that the variable $p$ appears in the energy function only via the quadratic form $\frac{1}{2}p^TM^{-1}p$. Thus, our modelling assumption allows us to approximate the mass matrix $M$ just from a set of trajectories, hence capturing the dependency of $H$ on the variable $p$ also outside $\Omega_{[0,T]}$.} Other possible improvements can be obtained when some symmetry structure is known for the Hamiltonian. On a similar direction, in Subsection \ref{se:noise}, we add some regularization based on other prior physical knowledge.

\begin{figure}
\centering
\begin{subfigure}{.45\textwidth}
    \centering
    \includegraphics[width=.98\textwidth]{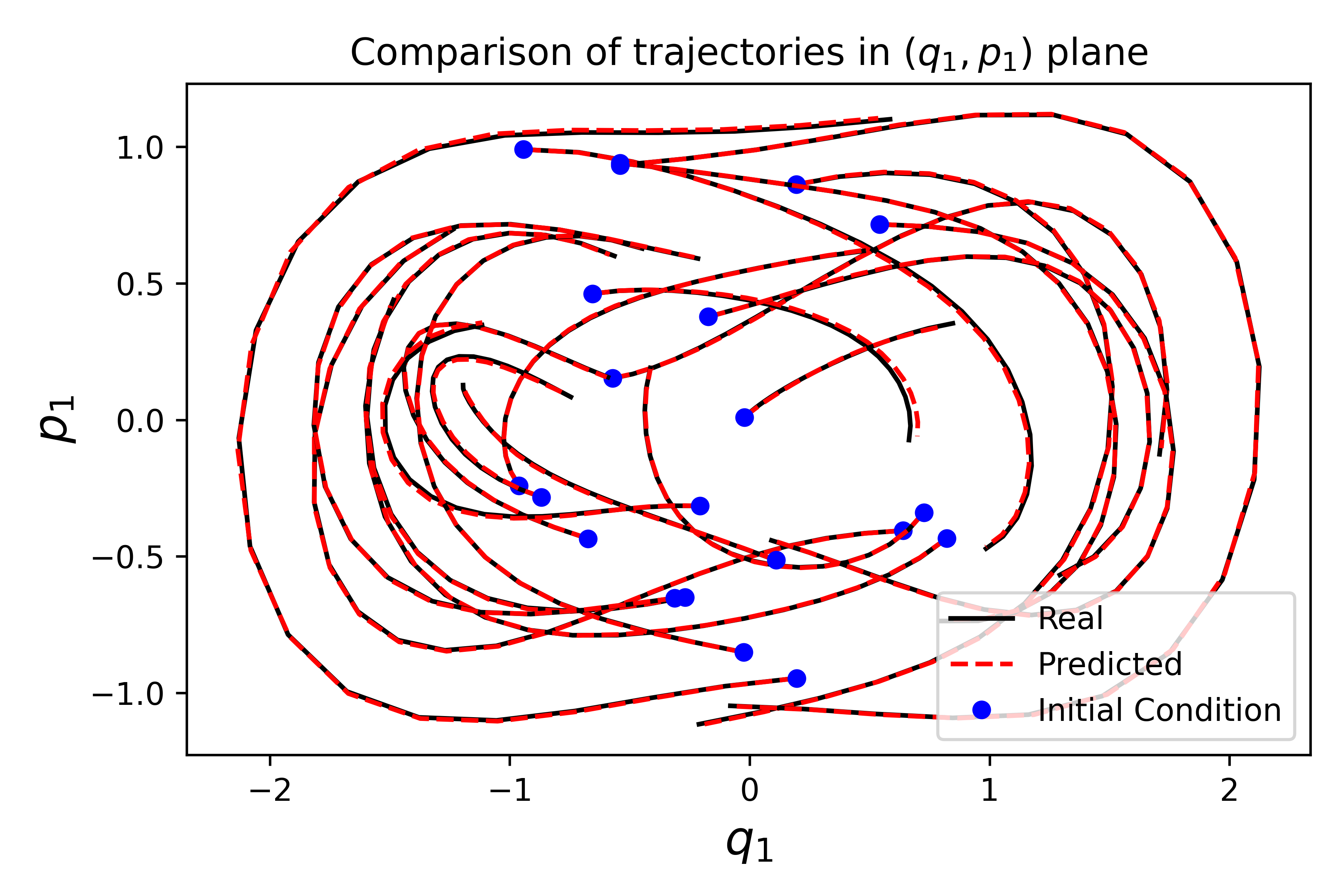}
    \caption{Projection on $(q_1,p_1)$}
\end{subfigure}
\begin{subfigure}{.45\textwidth}
    \centering
    \includegraphics[width=.98\textwidth]{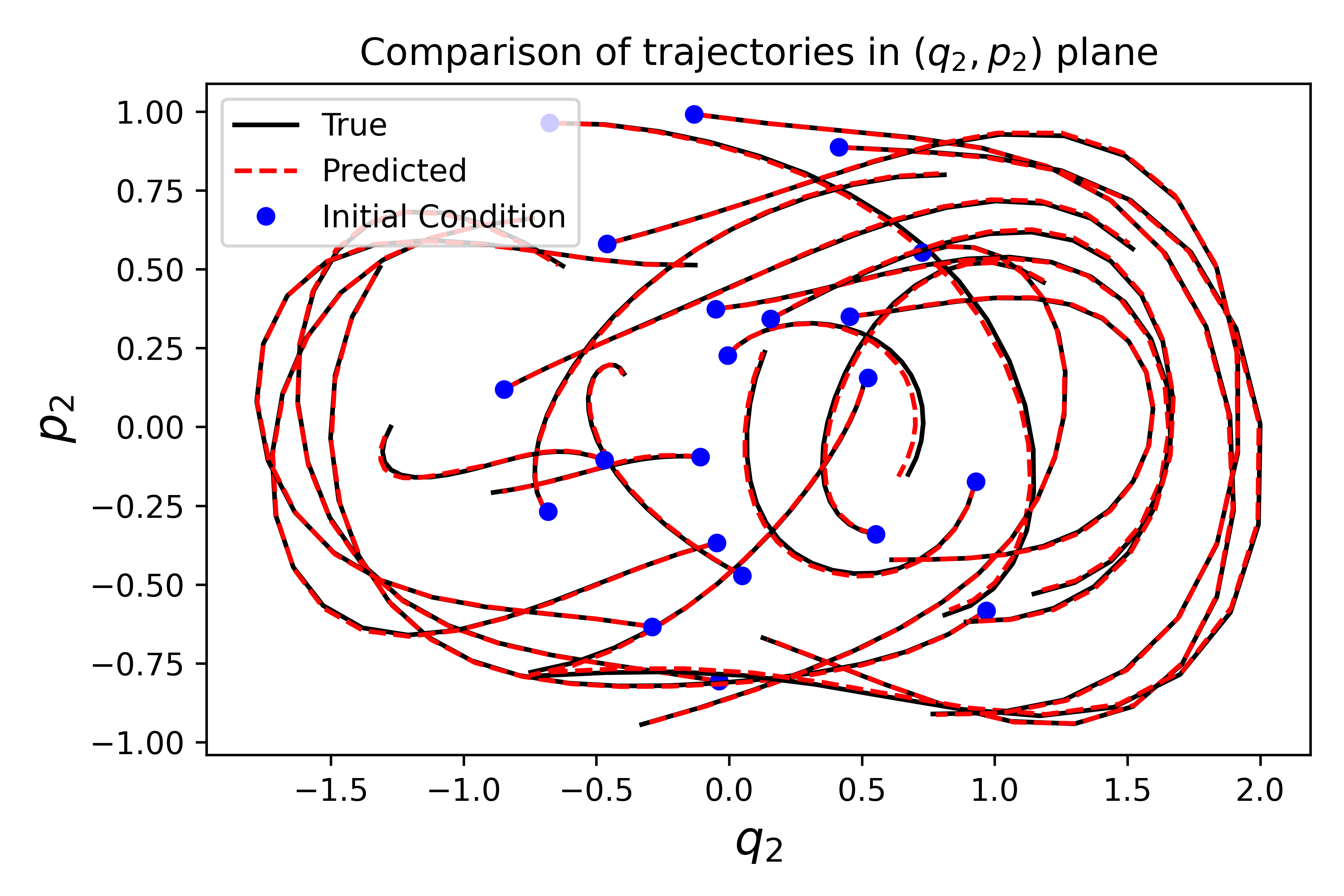}
    \caption{Projection on $(q_2,p_2)$}
\end{subfigure}
\caption{Comparison of real and predicted test trajectories for the Hamiltonian \eqref{eq:involved}. \red{In this case, for the potential energy we used a feedforward network with 3 hidden layers having respectively $100$, $50$ and $50$ neurons and $\tanh$ as activation function. The training integrator is Störmer-Verlet, with $M=6$ and final time $T=0.3$ and we use the Adam optimizer. The test trajectories, at $\tilde{M}=20$ uniformly distributed points in the time interval $[0,1]$, are obtained with ODE(5,4). These trajectories correspond to $\tilde{N}=100$ initial conditions on which the network has not been trained.}}
\label{fig:involved}
\end{figure}
We present in Figure~\ref{fig:involved} the comparison between ten learned trajectories and the corresponding exact ones of the Hamiltonian system $X_H\in\mathfrak{X}(\mathbb{R}^4)$ with Hamiltonian
\begin{equation}\label{eq:involved}
H(q,p) = \frac{1}{2}\begin{bmatrix} p_1 & p_2\end{bmatrix}^T\begin{bmatrix} 5 & -1 \\ -1 & 5 \end{bmatrix}\begin{bmatrix} p_1 \\ p_2\end{bmatrix} +  \frac{q_1^4+q_2^4}{4} + \frac{q_1^2+q_2^2}{2}.
\end{equation}
The training procedure of the network is based on 900 trajectories, sampled uniformly in 6 time instants, on the interval $0\leq t\leq 0.3$. \red{We remark that the training initial conditions are carefully chosen so that their associated trajectory segments well-capture the dynamics of interest.} \red{Figure \ref{fig:involved} collects test trajectories corresponding to the time interval $[0,1]$}. Since we are interested in approximating the Hamiltonian and not directly the trajectories, we are not constrained to evaluate the quality of the approximation with the same time integrator as the one used for training. In fact, these test trajectories have been generated with an embedded Runge--Kutta pair of order (5,4), with same relative and absolute accuracies for both the real and learned systems. Experimentally, it is clear that the qualitative behaviour of the Hamiltonian is well captured\red{, as we can see from Figure \ref{fig:involved}.} \red{To quantify the agreement of the prediction with the true Hamiltonian we report the $\mathcal{E}_1$ metric, as defined in \eqref{eq:MSETraj}, that is $6.59\cdot 10^{-5}$. Furthermore, the training loss is $4.62\cdot 10^{-7}$}.
\subsection{Robustness to noise and regularization}\label{se:noise}
In real world applications, data is contaminated by noise which usually comes from the measurement process. Thus, we need to test the robustness of the learning framework to the presence of noise in the training trajectories. To do so, we synthetically generate the trajectories as before, and then add random normal noise to all the points except the initial condition (for an averaging strategy that allows to deal even with perturbed initial conditions, see, e.g., \cite{chen2019symplectic}).
By construction, the network necessarily learns a Hamiltonian function, that is expected to generate trajectories close to the noisy ones. Since the training does not rely on clean trajectories, it is reasonable not to expect neither a loss value which is as small as in the absence of noise, nor a too accurate approximation of the Hamiltonian and the trajectories. Nevertheless, we aim for a learned Hamiltonian with level sets close to the exact ones, hence giving trajectories that resemble the true ones.
One way of improving the quality of the neural networks proposed here, is to make use of a priori known physical properties of the dynamical system. We use an approach based on soft constraints which means that we take the known physical properties into account by adding a regularization term in the cost function. An example of such a property could be one or more known conserved quantities, so called first integrals. Hamiltonian systems always have at least one first integral, namely the Hamiltonian function itself, but there might be additional independent ones. Enforcing the first integrals to be preserved or nearly preserved seems to be a reasonable strategy for obtaining improved qualitative behaviour of the resulting approximation as shown in the following example. 

\begin{figure}[h!]
\centering
\begin{subfigure}{.32\textwidth}
    \centering
    \includegraphics[width=\textwidth]{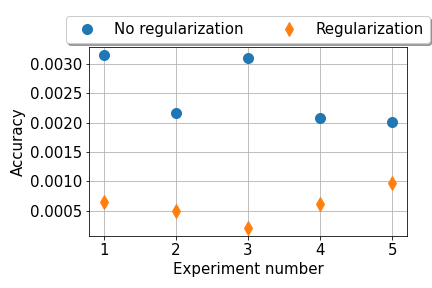}
    \caption{Case $\varepsilon=0.1$}
\end{subfigure}
\begin{subfigure}{.32\textwidth}
    \centering
    \includegraphics[width=\textwidth]{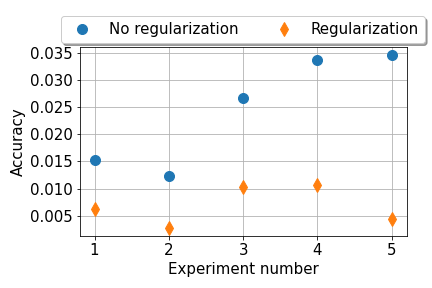}
    \caption{Case $\varepsilon=0.3$}
\end{subfigure}
\begin{subfigure}{.32\textwidth}
    \centering
    \includegraphics[width=\textwidth]{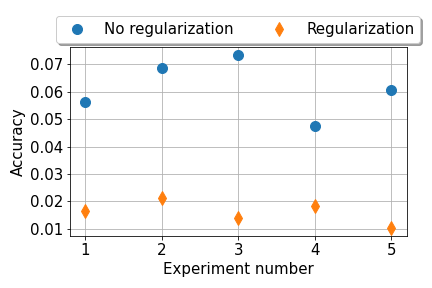}
    \caption{Case $\varepsilon=0.5$}
\end{subfigure}
\caption{5 repeated experiments for each perturbation regime. We plot on the $y$ axis the average accuracy, in terms of the $\mathcal{E}_1$ measure, obtained with the trained network, when compared with the real (non-noisy) trajectories.}
\label{fig:noisy}
\end{figure}

Consider a Hamiltonian system with Hamiltonian function $H:\mathbb{R}^{2n}\rightarrow\mathbb{R}$, and a functionally independent first integral $G$, i.e. $\nabla H(x)$ and $\nabla G(x)$ are never parallel. Consider the numerical integration $\hat{y}_k^j$, $j=1,...,M$ of the approximated Hamiltonian vector field $X_{H_{\Theta}}$,  starting at $\hat{y}_k^1 = x_k $. In the ideal case in which the learned Hamiltonian $H_{\Theta}$ coincides with $H$ and the numerical flow is replaced with the exact one, both $H$ and $G$ should be conserved. For this reason, we suggest adding to the loss function in equation \eqref{eq:loss} the following ``regularization" term:
\[
\mu\sum_{j\in \mathcal{I}}\left(G(\hat{y}_k^j)-G(x_k)\right)^2
\]
for all the training points $x_k$. Here $\mathcal{I}$ is a subset of indices contained in $\{1,..., M\}$, and $\mu$ is a regularization parameter that balances the importance of the preservation of the additional first integral against the perfect fitting of the training trajectories.
We test this regularization procedure with the Hamiltonian system $X_H\in\mathfrak{X}(\mathbb{R}^4)$ defined by
\[
H(q_1,q_2,p_1,p_2) = \frac{q_1^2+p_1^2}{2} + \frac{p_2^2}{2}+\frac{1}{2}q_2^2 + \frac{1}{4}q_2^4= h_1(q_1,p_1) + h_2(q_2,p_2).
\]
This system has $G(q,p):=h_1(q_1,p_1)$ as an additional independent first integral other than $H$.  We report in Figure \ref{fig:noisy} some plots of the obtained $\mathcal{E}_1$ values as defined in \eqref{eq:MSETraj}. In these experiments we add some random noise of the form $\varepsilon \delta$ to the points $y_i^j$ of the numerical trajectories, where $\delta\red{\sim}\mathcal{N}(0,1)$ \red{follows a standard normal distribution}. The same experiment is run $5$ times, and for each of these we plot the obtained $\mathcal{E}_1$ value. \red{For each experiment we generate new training and test trajectories, and these are used for both the regularized training and the non regularized one. Furthermore, each experiment has a different random initialization of the weights, which is however shared between the regularized and non regularized networks.}  We notice that with regularization we can consistently get a better error in terms of the $\mathcal{E}_1$ measure. There is not a huge difference between the results, \red{however. This suggests} that when prior information is known, it might be important to experiment with this kind of regularizing terms.
To conclude the Section, we highlight how remarkable it is that even without the regularization term, the trajectories are qualitatively well captured by the network and hence the test error is quite low. This is mostly due to the prior physical knowledge we impose on the learning procedure, i.e. that the vector field should be Hamiltonian. Indeed, since in the worst case the network approximates the wrong Hamiltonian, we always expect that it does not overfit the noisy trajectories, since they can not be learned exactly. On the other hand, without the prior knowledge of the Hamiltonian nature of the system, all the overfitting problems of standard neural networks reoccur and the risk of being closer to an interpolant of the noisy trajectories is higher.

\section{Learning constrained Hamiltonian systems} \label{se:constr}
The approximation of the Hamiltonians of constrained mechanical systems \red{with neural networks} has already been studied in the literature. Two main approaches can be identified. One of them is based on local coordinates on the constrained manifold (see, e.g., \cite{chen2019symplectic}, \cite{greydanus2019hamiltonian}) and the other uses ambient space coordinates and Lagrange multipliers (see \cite{finzi2020simplifying}). 
In principle, both the formulations apply to any constrained Hamiltonian system. However, as remarked in \cite{finzi2020simplifying}, the choice of a redundant system of coordinates usually gives a simpler expression for the Hamiltonian. This results in a more data efficient training procedure. In the second approach an embedded Runge--Kutta pair of order (5,4) is used to train the network. This choice inevitably leads to a drift from the constrained manifold during the training, even if it can be reduced by setting the tolerances of the integrator. However, in this way the cost of the integrator increases, hence this is not the most efficient way to preserve the constraints.

In this work, we use an alternative global formulation of the dynamics, as introduced in Section \ref{se:ham}. In principle this formulation adapts to any constrained Hamiltonian system whose configuration manifold is a submanifold of $\mathbb{R}^n$. Coupling this description of the dynamics with the learning framework introduced in Section \ref{se:hamNet}, their Hamiltonian functions can be approximated. To be more precise, one can use any numerical integrator to discretize the constrained trajectories and compare them with the training data. For example, Runge--Kutta 4 method can be used and this experimentally gives fast training procedures and accurate approximations of the Hamiltonian, as shown in the experiments of Subsection \ref{se:constrExp}.

We remark that in general numerical integrators do not preserve the geometry of the system and there might be a drift from the constrained manifold (see, e.g., \cite[Chapter 7]{hairer2006geometric}). Experimentally this does not seem to have a great impact on the quality of the predicted Hamiltonian in most of the cases. However, as we present in the numerical experiments with Lie group integrators, there might be situations in which one benefits from training the Hamiltonian with an integrator preserving the phase space. 
\red{Notice that the Hamiltonian that defines the dynamics has non-unique extension outside the phase space $\mathcal{M}=T^*\mathcal{Q}$. This is due to the projection matrix $P(q) = I_n - G(q)\left(G(q)^TG(q))\right)^{-1}G(q)^T$ appearing Equation \eqref{eq:ProjectionEq}, where $G(q)$ is the Jacobian matrix of the constraint function $g(q)$ defining $\mathcal{Q}$. This justifies investigating the importance of the preservation of the manifold $T^*Q$ in the training procedure.}


As introduced in Section \ref{se:ham}, in this work we assume that the constrained configuration manifold $\mathcal{Q}$ is known. Referring to equation (\ref{eq:ProjectionEq}), we notice that once the geometry is known, it is enough to specify the Hamiltonian function $H:T^*\mathcal{Q} \subseteq \mathbb{R}^{2n}\rightarrow\mathbb{R}$ in order to characterize the dynamics of a system. We show a setting in which the geometry can be preserved by Lie group integrators (see, \cite{iserles2000lie}, \cite{celledoni2014introduction}, \cite{celledoni2021}) focusing on the case $T^*\mathcal{Q}$ is homogeneous\footnote{A smooth manifold $\mathcal{M}$ is homogeneous if for any pair of points $m_1,m_2\in \mathcal{M}$ there is $g\in G$ such that $\psi(g,m_1)=m_2$, where $\psi:G\times\mathcal{M}\rightarrow \mathcal{M}$ is a Lie group action. In other words, $\psi$ is a transitive action.}. We see this even as an opportunity to study the behaviour of this class of methods in an applied framework and combined with neural networks. This geometric setup applies, for example, when $\mathcal{Q}$ is a homogeneous manifold and the transitive action $\psi:G\times \mathcal{Q}\rightarrow \mathcal{Q}$ defines, for any $q\in \mathcal{Q}$, a submersion $\psi_q:G\rightarrow \mathcal{Q}$ at the identity element $e\in G$ (see, e.g., \cite{brockett1972tangent}, \cite{celledoni2021dynamics}). Cartesian products of homogeneous manifolds are homogeneous too. Usually, multibody systems have constrained configuration manifolds given by cartesian products of $S^2$, $\mathbb{R}^k$, $SO(3)$ and $SE(3)$, which are respectively the special orthogonal and Euclidean groups. These are all homogeneous manifolds and so are their tangent and cotangent bundles.
\subsection{Lie group methods and neural networks}
Among the various classes of Lie group methods, we consider the Runge--Kutta--Munthe--Kaas (RKMK) methods and the commutator free ones (see, e.g., \cite{munthe-kaas99hor}, \cite{celledoni2003commutator}). The underlying idea of RKMK methods, applied to $F\in\mathfrak{X}(\mathcal{M})$\red{, with $\mathcal{M}$ an arbitrary homogeneous manifold,} is to express $F$ as $F\vert_m = \psi_*(f(m))\vert_m$. \red{Here} $\psi_*$ is the infinitesimal generator of $\psi$, a transitive Lie group action of $G$ on $\mathcal{M}$, and $f:\mathcal{M}\rightarrow\mathfrak{g}$ is a function that locally lifts the dynamics to the Lie algebra $\mathfrak{g}$ of $G$. On this linear space, we can perform a time step integration. We then map the result back to $\mathcal{M}$, and repeat this up to the final integration time.  More explicitly, let $\Delta t$ be the size of the uniform time step of the discretization, we then update $y_n\in \mathcal{M}$ to $y_{n+1}$ by
\begin{equation}\label{eq:update}
    \begin{cases}
    \gamma(0) = 0\in\mathfrak{g},\\
    \dot{\gamma}(t) = \text{dexp}_{\gamma(t)}^{-1}\circ f\circ \psi (\exp(\gamma(t)),y_n)\in T_{\gamma(t)}\mathfrak{g}, \\
    y_{n+1} = \psi(\exp(\gamma_1),y_n)\in \mathcal{M},
    \end{cases}
\end{equation}
where $\gamma_1\approx \gamma(\Delta t)\in\mathfrak{g}$ is computed with a Runge--Kutta method\red{, and $\text{dexp}^{-1}$ is the inverse of the differential of the exponential map $\text{exp}:\mathfrak{g}\rightarrow G$ as defined, for example, in \cite[Section 2.6]{iserles2000lie}.} We do not go into the details of commutator free methods, but the following development applies to them as well. In particular the function $f$ still plays a fundamental role.

We now present a natural way to combine the learning framework typical of unconstrained systems with Lie group integrators. This is done introducing a Lie group method during the learning procedure. Indeed, since we want to apply a Lie group integrator to deal with \red{nonlinear} geometries, we set $\Psi^{\Delta t}$, defined in equation (\ref{eq:psi}), to be the $\Delta t$ update given by some RKMK method. In other words, using the notation of equation (\ref{eq:update}), we get $\Psi^{\Delta t}(z) = \psi (\exp(\gamma_1),z)$ with $\gamma_1\in\mathfrak{g}$. 

The setting presented above for generic vector fields on homogeneous manifolds simplifies considerably in the presence of Hamiltonian systems. Indeed, for this type of systems, what is needed to fully determine the dynamics is the geometry given by $\mathcal{M}=T^*\mathcal{Q}$ and the scalar Hamiltonian function $H:\mathcal{M}\rightarrow \mathbb{R}$. In other words, we can think of the function $f:\mathcal{M}\rightarrow \mathfrak{g}$, that allows to express the vector field in terms of the infinitesimal generator of the action, as the result of an operator $F:C^{1}(\mathcal{M},\mathbb{R})\rightarrow  \{T^*\mathcal{M}\rightarrow \mathfrak{g}\}$ acting on a scalar function $H$. More explicitly, we can write $f=F[H]$ where $F$ and $H$ encode respectively the geometry and the dynamics of the system. This operator is not really necessary, but it clarifies considerably how the neural network comes into play in the learning framework. Indeed, because of this construction, we can write the numerical flow $\Psi^{\Delta t}$ as the map sending $y_n$ into $y_{n+1}=\psi(\exp(\gamma_{\Delta t,y_n}),y_n)$ with $\gamma_{\Delta t,y_n}$ being an approximation of the solution $\gamma(\Delta t)$ of the following initial value problem \[
\begin{cases}
\dot{\gamma}(t) = \mathrm{dexp}_{\gamma(t)}^{-1}\circ F[H_{\Theta}]\circ \psi (\exp(\gamma(t)),y_n)\in T_{\gamma(t)}\mathfrak{g},\\
\gamma(0)=0\in\mathfrak{g}.
\end{cases}
\]
Here $H_{\Theta}$ is the approximation of the Hamiltonian given by the current weights $\Theta$ of the neural network. Thus, applying a particular family of geometric numerical integrators, we can directly study some constrained systems with the same ideas coming from learning unconstrained ones. Since following this procedure the geometry is preserved, one can \red{consider replacing} the Euclidean distance in the loss function defined in equation (\ref{eq:loss}) with a Riemannian metric of the constrained manifold. \red{This would bring to distances between points that correspond to the length of the minimal geodesic connecting them, which is in general different from the length of the segment in the ambient space having them as extrema.} In the remaining part of the Section, we specialize this reasoning to mechanical systems defined on copies of $T^*S^2$. We focus on a chain of spherical pendula, but the geometric setting applies also to other systems (see, e.g., \cite[Section 10.5]{lee2017global}).
\subsection{Mechanical systems on $(T^*S^2)^k$}\label{se:spherical}
As anticipated in the introductory Section, in this geometric setting we are not involving symplectic integrators and we do not assume to have a separable Hamiltonian anymore. Thus, we now model a more general family of Hamiltonians as
\begin{equation}\label{ham}
   H(q,p) = \frac{1}{2}p^TM^{-1}(q)p + V(q).
\end{equation}
We model the potential energy as before, however we need an alternative strategy for the inverse of the mass matrix, which is no longer assumed to be constant. Based on the problem, one can choose various parametrizations of the mass matrix or its inverse. We decide to specialize the architecture based on the fact that the geometry of the system is known to be $\mathcal{M}=(T^*S^2)^k$, where $S^2\subset\mathbb{R}^3$. 
\red{We coordinatize $\mathcal{M}$ with $(q,p)=(q_1,\dots,q_k,p_1,\dots,p_k)\in \mathbb{R}^{6k}$.}
In this case, when \red{$p\in \mathbb{R}^{3k}$} is intended as the vector of linear momenta, the matrix $M(q)$ in equation (\ref{ham}) is a block matrix, with
\[
i,j = 1,...,k,\quad \mathbb{R}^{3\times 3}\ni M(q)_{ij} = \begin{cases} m_{ii}I_3,\quad i=j\\
m_{ij}(I_3-q_iq_i^T),\quad \text{otherwise,}
\end{cases}
\]
see \cite[Section 8.3.3]{lee2017global} for further details. Here, the matrix having constant entries $m_{ij}$ is symmetric and positive definite. For this reason, we leverage this form of the kinetic energy and learn a constant matrix $A\in\mathbb{R}^{k\times k}$ and a vector $b\in\mathbb{R}^k$ so that
\begin{equation}\label{eq:massmat}
\begin{bmatrix}
m_{11} & ... & m_{1k}\\
m_{21} & ... & m_{2k}\\
\vdots & \vdots & \vdots \\
m_{k1} & ... & m_{kk}
\end{bmatrix} \approx A^TA + 
\begin{bmatrix}
\tilde{b}_{1} & 0 & ... & 0 \\
0 & \tilde{b}_2 & \ddots & \vdots \\
\vdots & \ddots & \ddots & 0 \\
0 & ... & 0 & \tilde{b}_k
\end{bmatrix}
\end{equation}
where $\tilde{b}_i := \max{(0,b_i)}$ are terms added to promote the positive definiteness of the right-hand side. We tested also elevating to the second power the $b_i$ instead of taking the maximum with $0$, but we got better results with the choice presented in equation \eqref{eq:massmat}. The matrix on the left-hand side of equation \eqref{eq:massmat} is exactly the one appearing in Hamiltonian formulations with Cartesian coordinates, as the one used in \cite{finzi2020simplifying}.

For the spherical pendulum we have $k=1$ and hence the Hamiltonian dynamics is defined on its cotangent bundle $T^*S^2$, which is a homogeneous manifold. This can be obtained thanks to the transitivity of the group action 
\[
\Psi : SE(3)\times T^*S^2\rightarrow T^*S^2,\;\;((R,r),(q,p^T))\mapsto (Rq,(Rp+r\times Rq)^T),
\]
where the transpose comes from the usual interpretation of covectors as row vectors. \red{As in \cite[Chapter 6]{holm2011geometric}, we represent a generic element of the special Euclidean group $G=SE(3)$ as an ordered pair $(R,r)$, where $R\in SO(3)$ is a rotation matrix and $r\in\mathbb{R}^3$ is a vector.} With this specific choice of the geometry, the formulation presented in equation (\ref{eq:ProjectionEq}) simplifies considerably. Indeed $P(q) = I_3-qq^T$ which implies $W(q,p) = pq^T-qp^T$. Replacing these expressions in (\ref{eq:ProjectionEq}) and using the triple product rule we end up with the following set of ODEs
\begin{equation}\label{eq:odepend}
\begin{cases}
\dot{q} &= (I-qq^T)\partial_pH(q,p)\\
\dot{p} &= -(I-qq^T)\partial_qH(q,p) + \partial_pH(q,p)\times (p\times q).
\end{cases}
\end{equation}
This vector field $X(q,p)$ can be expressed as $\psi_*(F[H](q,p))(q,p)$ with
\[
\psi_*((\xi,\eta))(q,p) = (\xi\times q,\xi\times p + \eta\times q), \red{\quad (\xi,\eta)\in \mathfrak{g= se}(3)}
\]
and
\[
F[H](q,p) = (\xi,\eta)=\left(q\times \frac{\partial H(q,p)}{\partial p},\frac{\partial H(q,p)}{\partial q}\times q + \frac{\partial H(q,p)}{\partial p}\times p \right).
\]

A similar reasoning can be extended to a chain of $k$ connected pendula, and hence to a system on $(T^*S^2)^k$. The main idea is to replicate both the equations \eqref{eq:odepend} and the expression $F[H]$ for all the $k$ copies of $T^*S^2$. A more detailed explanation can be found in \cite{celledoni2021}.

\begin{figure}[h!]
\centering
\begin{subfigure}{0.49\linewidth}
\centering
\includegraphics[width=\textwidth]{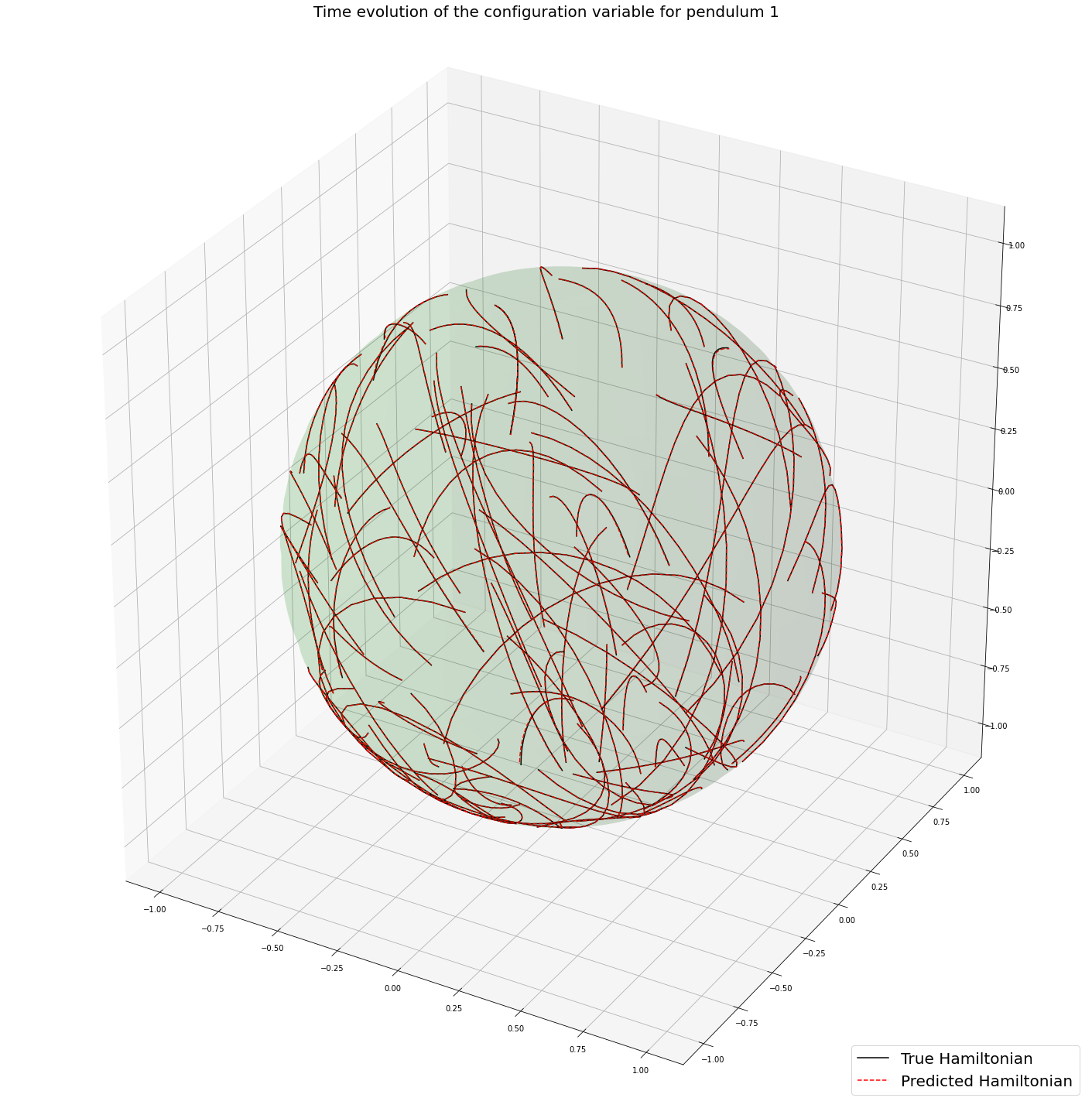}
\end{subfigure}
\begin{subfigure}{0.49\linewidth}
\centering
\includegraphics[width=\textwidth]{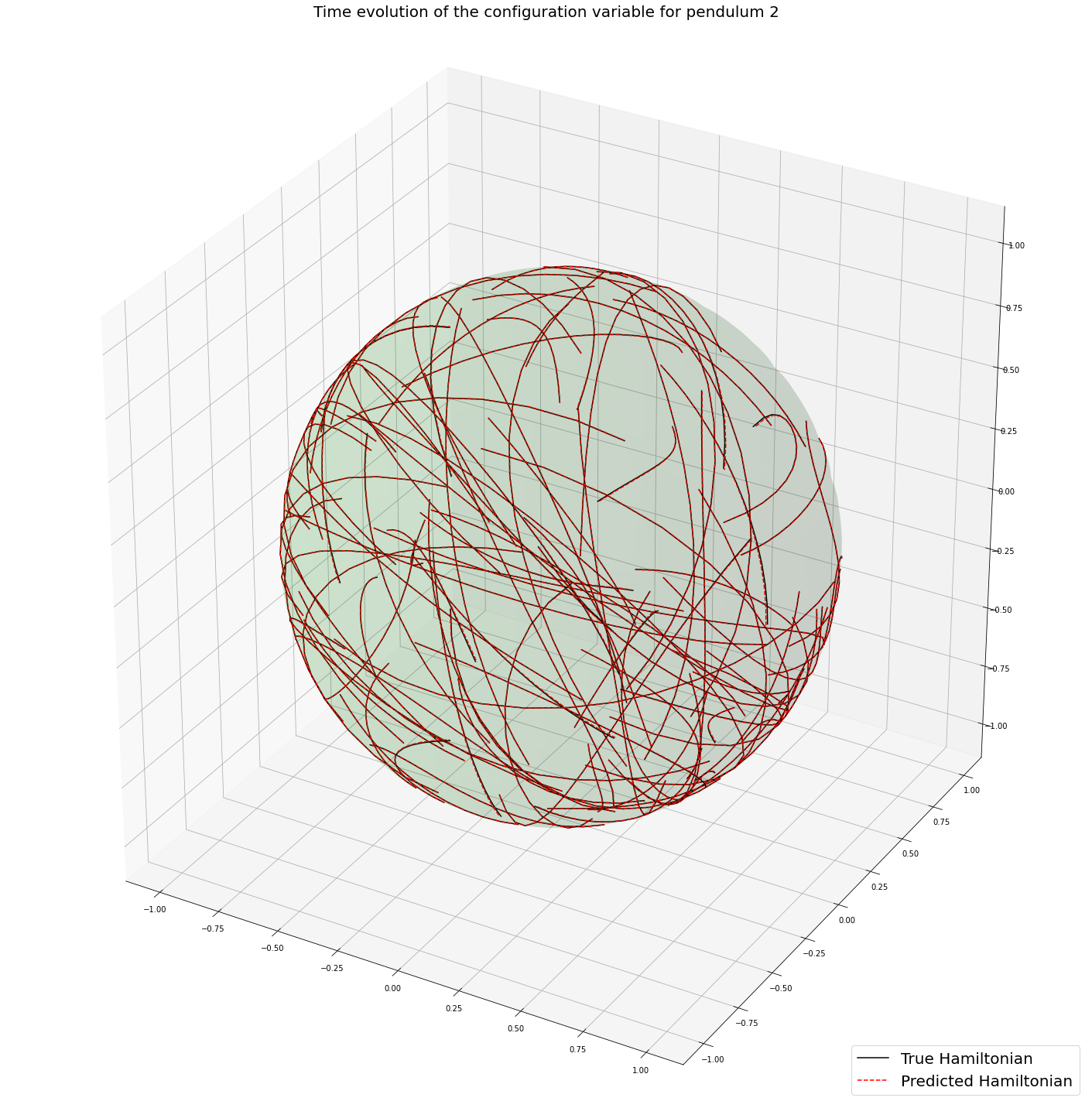}
\end{subfigure}
\caption{Comparison between 100 test trajectories obtained with the true Hamiltonian $H$ and the predicted one $H_{\Theta}$. To train $H_{\Theta}$, a Lie group method is used. This gives $\mathcal{E}_1=2.65\cdot 10^{-6}$ and a final training loss of $1.6
\cdot 10^{-9}$.}
\label{fig:2pend}
\end{figure}

We present in Figure \ref{fig:2pend} the results obtained for the training of a double pendulum, i.e. $k=2$. To train the network, we generate a set of $N=500$ training trajectories with the embedded Runge--Kutta pair of order (5,4) of $\texttt{SciPy}$. The final integration time is $T=0.1$ and $M=5$. \red{To model the potential energy, we use a feedforward network with 3 hidden layers of 100 neurons each.} In the plots we show the configuration variables, $q_1,\,q_2\in S^2$, obtained for 100 test trajectories in the time interval $[0,1]$, where the network $H_{\Theta}$ has been trained with a commutator free method of order 4.

\subsection{Experimental study of the learning procedure}\label{se:constrExp}
We investigate the influence of the training setup on the \red{error measures $\mathcal{E}_1$, $\mathcal{E}_2$, defined in \eqref{eq:avgError}, and on the training loss}. More precisely, we test how the parameters $M$, $N$, \red{the noise magnitude} and the training integrator affect the performance of the network. \red{We quantify the magnitude of noise in the training trajectories with a parameter $\varepsilon>0$, as in Subsection \ref{se:noise}.} The integrators that we study are Lie Euler, explicit Euler (both of order 1), commutator free and Runge--Kutta (both of order 4). In particular, Lie Euler and commutator free methods preserve the phase space $\mathcal{M}$ up to machine accuracy. To get a sufficient sample of experiments, we repeat all the tests 5 times, and look at the medians and geometric means\red{\footnote{\red{The choice of geometric means is because of the exponential nature of the error measures and the training loss.}}} of the obtained results.  \red{To be precise, we test $N\in\{50,500,1000,1500\}$, $M\in\{2,3,5\}$, and $\varepsilon\in\{0,0.001,0.01,0.1\}$. Therefore, we perform a total of $960$ experiments, and also here the potential energy is modelled with a feedforward network of 3 hidden layers having 100 neurons each. Furthermore, for the four experiments performed varying just the integrator, and with the other parameters fixed, the network's weights are initialized to be the same, and also the training and test initial conditions are the same.} For all these experiments, we focus on the single spherical pendulum, we keep the final training time to $T=0.1$, and we don't use regularization terms. The training trajectories have been generated with the $\texttt{SciPy}$ implementation of the Dormand-Prince pair of order (5,4) with strict tolerance.
\begin{table}[h!]
\red{
\centering
\small
\begin{tabular}{|c|c|c|c|c|} 
\hline
Order & Integrator & $\mathcal{E}_1$ & $\mathcal{E}_2$ & Training Loss \\
\hline
\hline 1 & EE & 5.7e-5 & 1.13e-2 & 2.12e-6 \\
\hline
1 & LE & 4.9e-5 & 1.07e-2 & 1.17e-6 \\
\hline
4 & RK4 & 1.12e-5 & 3.83e-3 & 2.63e-7 \\
\hline
4 & CF4 & 1.12e-5 & 3.85e-3 & 2.64e-7\\
\hline
\end{tabular}
\caption{\red{In this table we report the geometric means of the quantities $\mathcal{E}_1$, $\mathcal{E}_2$ and the training loss. Here we average over all the 240 experiments that have the same integrator. We denote the four integrators with EE (explicit Euler), LE (Lie Euler), RK4 (Runge--Kutta 4), and CF4 (commutator free 4).}}
\label{ta:order}}
\end{table}
\begin{figure}[h!]
    \centering
    \includegraphics[trim={0.5 1.5cm 0.5 1.5cm},clip,width=1.0\textwidth]{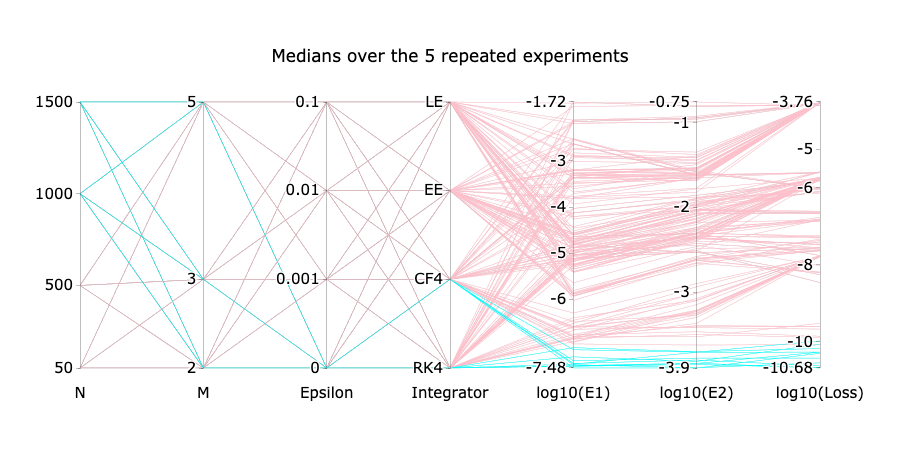}
    \caption{\red{This is a parallel coordinate plot reporting the dependencies of $\mathcal{E}_1$, $\mathcal{E}_2$ and the training loss on the parameters $N$, $M$, $\varepsilon$ and on the integrator. Each coloured polyline corresponds to the median over 5 experiments, i.e. same $N$, $M$, $\varepsilon$ and same integrator. The lines in cyan color represent the combinations giving $\mathcal{E}_1<10^{-7}$.}}
    \label{fig:parcoords}
\end{figure}
\red{As shown in Table \ref{ta:order}, the order of the numerical integrator used to train the network plays an important role. Indeed, we get results that are similar for methods of the same order, but there is a noticeable decay in the errors and in the loss when we increase the order from one to four. As highlighted in \cite{zhu2020deep}, this effect can be explained with a standard argument of backward error analysis, see e.g. \cite[Chapter 9]{hairer2006geometric}. From the results reported in Table \ref{ta:order} we see that the local error of the integrator is more important than the preservation of the geometry. Therefore, even if from a theoretical point of view it seems relevant to remain on the manifold during the training, in practice this does not seem to be very important in the particular experiment considered here. In Figure \ref{fig:parcoords}, we plot the dependencies of $\mathcal{E}_1$, $\mathcal{E}_2$ and the training loss, on $N$, $M$, $\varepsilon$ and the integrator. We notice that values of $\mathcal{E}_1$ below a threshold of $10^{-7}$ can be reached only with integrators of order four and with the smallest value of $\varepsilon$. The interplay of $N$, $M$ and $\varepsilon$ is further investigated in Table \ref{ta:comparison}. An interactive version of Figure \ref{fig:parcoords}, together with other parallel coordinate plots, can be found at the GitHub Page \url{https://davidemurari.github.io/learningConstrainedHamiltonians/}, while the dataset is available in the GitHub repository associated to the paper.}

\red{We conclude this parameter study considering separately the case with and without noise, $\varepsilon> 0$ and $\varepsilon=0$ respectively. The results are reported in Table \ref{ta:comparison}. In general the lowest values of $\mathcal{E}_1$ are obtained with high $N$. For the model under consideration, $N=1000$ seems already high enough to achieve good results. Regarding $M$, Table \ref{ta:comparison} shows that to achieve lower values of $\mathcal{E}_1$ in the presence of noise, one needs to adopt a higher $M$. On the other hand, in the absence of noise it seems important to have a high $M$ only for low order integrators. Finally, as may be expected, even if this Table does not distinguish among the different magnitudes of the noise, we see that with $\varepsilon=0$ better results can be achieved. }

\begin{table}[]
\centering
\red{
\begin{tabular}{|cccccccccc|}
\hline
\multicolumn{10}{|c|}{Without noise}     \\ \hline \hline
\multicolumn{5}{|c||}{Integrator of order 1}                                        & \multicolumn{5}{c|}{Integrator of order 4}                                 \\ \hline
\multicolumn{1}{|c|}{$N$}  & \multicolumn{1}{c|}{$M$} & \multicolumn{1}{c|}{Int.} & \multicolumn{1}{c|}{$\mathcal{E}_1$} & \multicolumn{1}{c||}{$\mathcal{E}_2$} & \multicolumn{1}{c|}{$N$}  & \multicolumn{1}{c|}{$M$} & \multicolumn{1}{c|}{Int.} & \multicolumn{1}{c|}{$\mathcal{E}_1$} & $\mathcal{E}_2$ \\ \hline
\multicolumn{1}{|c|}{1500} & \multicolumn{1}{c|}{5}   & \multicolumn{1}{c|}{LE}   & \multicolumn{1}{c|}{1e-6}            & \multicolumn{1}{c||}{2.4e-3}          & \multicolumn{1}{c|}{1500} & \multicolumn{1}{c|}{2}   & \multicolumn{1}{c|}{CF4}  & \multicolumn{1}{c|}{3.2e-8}          & 1.3e-4          \\ \hline
\multicolumn{1}{|c|}{1000} & \multicolumn{1}{c|}{5}   & \multicolumn{1}{c|}{LE}   & \multicolumn{1}{c|}{1e-6}            & \multicolumn{1}{c||}{2.3e-3}          & \multicolumn{1}{c|}{1500} & \multicolumn{1}{c|}{5}   & \multicolumn{1}{c|}{RK4}  & \multicolumn{1}{c|}{3.3e-8}          & 1.4e-4          \\ \hline
\multicolumn{1}{|c|}{1000} & \multicolumn{1}{c|}{5}   & \multicolumn{1}{c|}{EE}   & \multicolumn{1}{c|}{1e-6}            & \multicolumn{1}{c||}{2.4e-3}          & \multicolumn{1}{c|}{1500} & \multicolumn{1}{c|}{3}   & \multicolumn{1}{c|}{RK4}  & \multicolumn{1}{c|}{3.4e-8}          & 1.4e-4          \\ \hline
\multicolumn{1}{|c|}{1500} & \multicolumn{1}{c|}{5}   & \multicolumn{1}{c|}{EE}   & \multicolumn{1}{c|}{1e-6}            & \multicolumn{1}{c||}{2.5e-3}          & \multicolumn{1}{c|}{1500} & \multicolumn{1}{c|}{2}   & \multicolumn{1}{c|}{RK4}  & \multicolumn{1}{c|}{3.6e-8}          & 1.5e-4          \\ \hline
\multicolumn{1}{|c|}{500}  & \multicolumn{1}{c|}{5}   & \multicolumn{1}{c|}{LE}   & \multicolumn{1}{c|}{2e-6}            & \multicolumn{1}{c||}{2.6e-3}          & \multicolumn{1}{c|}{1500} & \multicolumn{1}{c|}{3}   & \multicolumn{1}{c|}{CF4}  & \multicolumn{1}{c|}{3.7e-8}          & 1.5e-4          \\ \hline \hline
\multicolumn{10}{|c|}{With noise}        \\ \hline \hline
\multicolumn{5}{|c||}{Integrator of order 1}                                       & \multicolumn{5}{c|}{Integrator of order 4}                                       \\ \hline
\multicolumn{1}{|c|}{$N$}  & \multicolumn{1}{c|}{$M$} & \multicolumn{1}{c|}{Int.} & \multicolumn{1}{c|}{$\mathcal{E}_1$} & \multicolumn{1}{c||}{$\mathcal{E}_2$} & \multicolumn{1}{c|}{$N$}  & \multicolumn{1}{c|}{$M$} & \multicolumn{1}{c|}{Int.} & \multicolumn{1}{c|}{$\mathcal{E}_1$} & $\mathcal{E}_2$ \\ \hline
\multicolumn{1}{|c|}{1500} & \multicolumn{1}{c|}{5}   & \multicolumn{1}{c|}{LE}   & \multicolumn{1}{c|}{1.2e-5}          & \multicolumn{1}{c||}{6.6e-3}          & \multicolumn{1}{c|}{1500} & \multicolumn{1}{c|}{5}   & \multicolumn{1}{c|}{RK4}  & \multicolumn{1}{c|}{5e-6}            & 3.4e-3          \\ \hline
\multicolumn{1}{|c|}{1500} & \multicolumn{1}{c|}{5}   & \multicolumn{1}{c|}{EE}   & \multicolumn{1}{c|}{1.2e-5}          & \multicolumn{1}{c||}{6.3e-3}          & \multicolumn{1}{c|}{1500} & \multicolumn{1}{c|}{5}   & \multicolumn{1}{c|}{CF4}  & \multicolumn{1}{c|}{6e-6}            & 4.1e-3          \\ \hline
\multicolumn{1}{|c|}{1000} & \multicolumn{1}{c|}{5}   & \multicolumn{1}{c|}{LE}   & \multicolumn{1}{c|}{1.5e-5}          & \multicolumn{1}{c||}{6.6e-3}          & \multicolumn{1}{c|}{1000} & \multicolumn{1}{c|}{5}   & \multicolumn{1}{c|}{CF4}  & \multicolumn{1}{c|}{7e-6}            & 3.8e-3          \\ \hline
\multicolumn{1}{|c|}{1000} & \multicolumn{1}{c|}{5}   & \multicolumn{1}{c|}{EE}   & \multicolumn{1}{c|}{1.6e-5}          & \multicolumn{1}{c||}{6.8e-3}          & \multicolumn{1}{c|}{1000} & \multicolumn{1}{c|}{5}   & \multicolumn{1}{c|}{RK4}  & \multicolumn{1}{c|}{8e-6}            & 4.3e-3          \\ \hline
\multicolumn{1}{|c|}{500}  & \multicolumn{1}{c|}{5}   & \multicolumn{1}{c|}{LE}   & \multicolumn{1}{c|}{1.8e-5}          & \multicolumn{1}{c||}{6.7e-3}          & \multicolumn{1}{c|}{1000} & \multicolumn{1}{c|}{3}   & \multicolumn{1}{c|}{RK4}  & \multicolumn{1}{c|}{8e-6}            & 4.2e-3          \\ \hline
\end{tabular}
\caption{\red{In this Table we report the combinations that give the 5 best values of $\mathcal{E}_1$, together with the corresponding value  $\mathcal{E}_2$. These are the geometric means among all the experiments. The two tables compare the performance on data with and without noise.}}
\label{ta:comparison}}
\end{table}

We also point out that the experiments were performed for short integration times, where not only symplectic integrators can generate physically meaningful trajectories. It would be interesting to explore the performance of sympelctic and constraint preserving integrators in this setting (see, e.g., \cite{andersen1983rattle}) and we defer this to further work.

\red{Besides the theoretical aspect of the non-uniqueness of the extension of the dynamics outisde of $\mathcal{M}\subset\mathbb{R}^{2n}$, we now report a numerical experiment where the preservation of the geometry during the training is beneficial.} We consider again a simple spherical pendulum and we assume to know that the potential energy is linear. We hence impose this prior information on the architecture of the network. Due to the problem's simplicity, we aim to reach very low $\mathcal{E}_1$ and $\mathcal{E}_2$ values. Training the same architecture for 200 epochs, both with Runge--Kutta and commutator free methods of order 4, we get the results in Table \ref{tab:comparison}. Indeed the geometric integrator outperforms the classical Runge-Kutta method in this experiment. 
\begin{table}[h!]
    \centering
    \begin{tabular}{|c|c|c|}
    \hline
        Numerical method in the training  & $\mathcal{E}_1$ & $\mathcal{E}_2$ \\
        \hline \hline
        Runge-Kutta of order 4 & 4.2e-12 & 1.5e-6\\
        \hline
        Commutator free of order 4 &  1.1e-14 &2.5e-7\\
        \hline
    \end{tabular}
    \caption{Comparison of the accuracy measures $\mathcal{E}_1$ and $\mathcal{E}_2$ obtained with the two integrators. \red{These results are obtained imposing the linear structure of the potential energy on the network modelling the Hamiltonian of the spherical pendulum. The kinetic energy has been modelled as in previous experiments.}}
    \label{tab:comparison}
\end{table}

This experiment suggests that the choice of an integrator that does not fully exploit the available information, like the geometry, \red{might limit the quality of the obtained approximations. For those cases in which one is interested in as accurate as possible predictions, this might be a relevant issue.} 

The experiments performed lead to the conclusion that modelling multi-body systems with neural networks can be a valuable approach. However, to better leverage the approximation capabilities of machine learning techniques (see, e.g., \cite{hornik1991approximation}, \cite{cybenko1989approximation}) we believe that a deeper investigation and understanding of how they interface with physical models is necessary.

\section*{Disclosure statement}
No potential conflict of interest was reported by the author(s).
\section*{Acknowledgements}

\begin{wrapfigure}{r}{0.15\textwidth}
    \centering
    \vspace{-10pt}
    \includegraphics[width=0.12\textwidth]{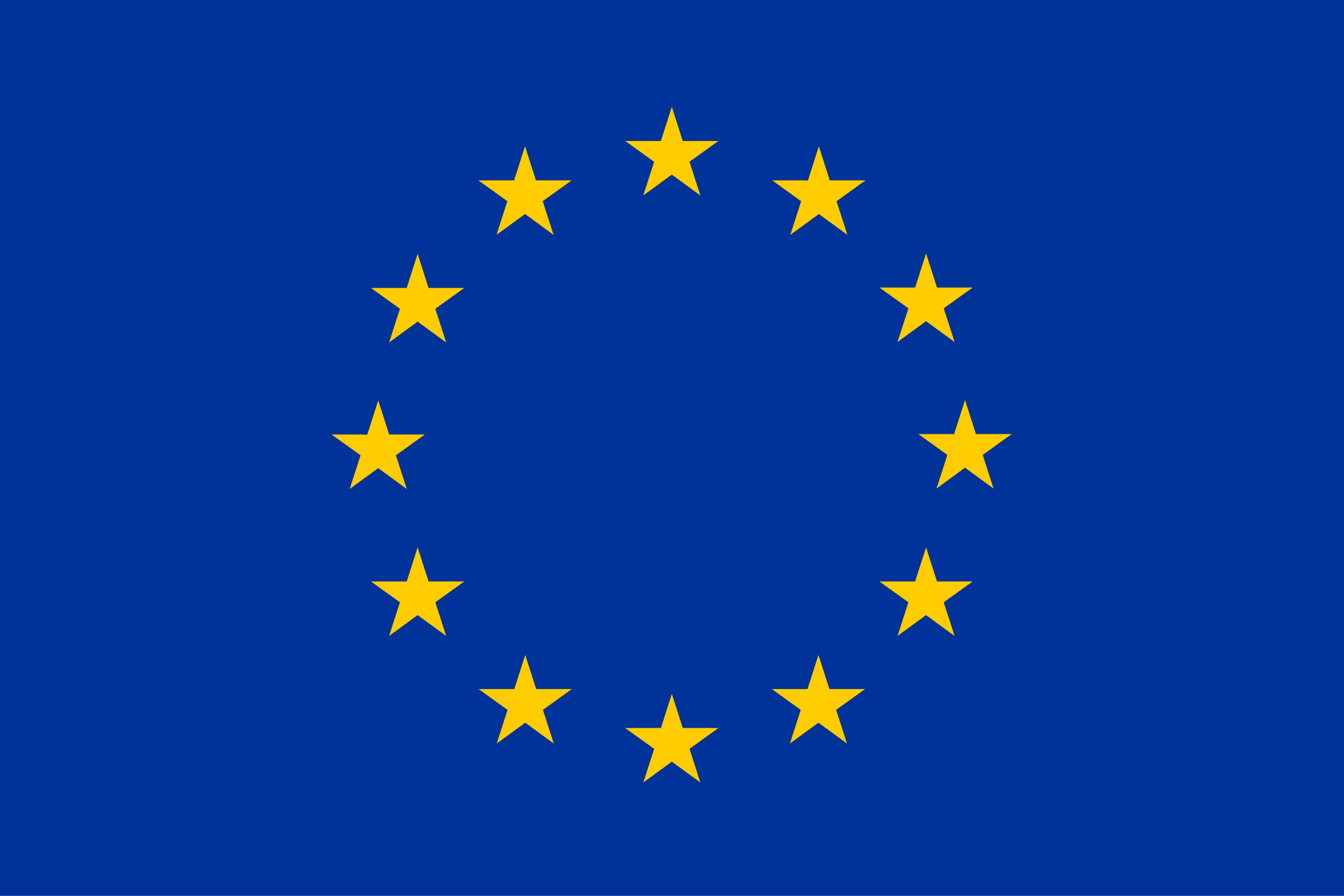}
\end{wrapfigure}
This project has received funding from the European Union’s Horizon 2020 research and innovation programme under the Marie Skłodowska-Curie grant agreement No 860124. \newline\newline 
The authors are grateful to Ergys Çokaj for the valuable discussions in the early stage of this work. The authors would like to thank the Isaac Newton Institute for Mathematical Sciences, Cambridge, for support and hospitality during the programme "Mathematics of Deep Learning".

\bibliographystyle{elsarticle-harv}
\bibliography{main}


\end{document}